\documentclass[12pt,leqno]{amsart}
\usepackage{amssymb,verbatim,enumerate,ifthen,cite}
\usepackage[mathscr]{eucal}
\usepackage[utf8]{inputenc}
\usepackage[T1]{fontenc}
\def\N{\mathbb{N}}
\def\R{\mathbb{R}}
\def\Q{\mathbb{Q}}

\def\K{\mathbb{K}}

\newtheorem{theorem}{Theorem}[section]
\newtheorem*{theorem*}{Theorem}
\def\Thm#1#2{\ifthenelse{\equal{#1}{*}}{\begin{theorem*}#2\end{theorem*}}
	{\begin{theorem}\label{T#1}#2\end{theorem}}}

\def\thm#1{Theorem~\ref{T#1}}

\newtheorem{proposition}[theorem]{Proposition}
\newtheorem*{proposition*}{Proposition}
\def\Prp#1#2{\ifthenelse{\equal{#1}{*}}{\begin{proposition*}#2\end{proposition*}}
	{\begin{proposition}\label{P#1}#2\end{proposition}}}

\newtheorem{corollary}[theorem]{Corollary}
\newtheorem*{corollary*}{Corollary}
\def\Cor#1#2{\ifthenelse{\equal{#1}{*}}{\begin{corollary*}#2\end{corollary*}}
	{\begin{corollary}\label{C#1}#2\end{corollary}}}

\newtheorem{lemma}[theorem]{Lemma}
\newtheorem*{lemma*}{Lemma}
\def\Lem#1#2{\ifthenelse{\equal{#1}{*}}{\begin{lemma*}#2\end{lemma*}}
	{\begin{lemma}\label{L#1}#2\end{lemma}}}
\def\lem#1{Lemma~\ref{L#1}}

\theoremstyle{definition}
\newtheorem{remark}[theorem]{Remark}
\newtheorem*{remark*}{Remark}
\def\Rem#1#2{\ifthenelse{\equal{#1}{*}}{\begin{remark}\rm #2\end{remark}}
	{\begin{remark}\label{R#1}\rm #2\end{remark}}}

\newtheorem{statement}[theorem]{Statement}
\newtheorem*{statement*}{Statement}
\def\Stm#1#2{\ifthenelse{\equal{#1}{*}}{\begin{statement*}#2\end{statement*}}
	{\begin{statement}\label{s#1}#2\end{statement}}}

\newtheorem{example}[theorem]{Example}
\newtheorem*{example*}{Example}
\def\Exa#1#2{\ifthenelse{\equal{#1}{*}}{\begin{example*}\rm #2\end{example*}}
	{\begin{example}\label{Ex#1}\rm #2\end{example}}}

\def\eq#1{{\rm(\ref{E#1})}}
\def\Eq#1#2{\ifthenelse{\equal{#1}{*}}
	{\begin{equation*}\begin{aligned}[]#2\end{aligned}\end{equation*}}
	{\begin{equation}\begin{aligned}[]\label{E#1}#2\end{aligned}\end{equation}}}

\def\comment#1{}

\begin{document}
	\begin{flushright}
	\end{flushright}
	\vspace{5mm}
	
	\date{\today}
	
	\title{On convexity properties with respect to a Chebyshev system}
	
	\author[Zs. P\'ales]{Zsolt P\'ales}
	\address[Zs. P\'ales]{Institute of Mathematics, University of Debrecen, Hungary}
	\email{pales@science.unideb.hu}
	\author[M. K. Shihab]{Mahmood Kamil Shihab}
	\address[M. K. Shihab]{Doctoral School of Mathematical and Computational Sciences, University of Debrecen, Hungary; Department of Mathematics, College of Education for Pure Sciences, University of Kirkuk, Iraq}
	\email{mahmood.kamil@science.unideb.hu; mahmoodkamil30@uokirkuk.edu.iq}
	
	\thanks{The research of the first author was supported by the K-134191 NKFIH Grant.}

\begin{abstract}
The main purpose of this paper is to introduce various convexity concepts in terms of a positive Chebyshev system $\omega$ and give a systematic investigation of the relations among them. We generalize a celebrated theorem of Bernstein--Doetsch to the setting of $\omega$-Jensen convexity. We also give sufficient conditions for the existence of discontinuous $\omega$-Jensen affine functions. The concept of Wright convexity is extended to the setting of Chebyshev systems, as well, and it turns out to be an intermediate convexity property between $\omega$-convexity and $\omega$-Jensen convexity. For certain Chebyshev systems, we generalize the decomposition theorems of Wright convex and higher-order Wright convex functions obtained by C.\ T.\ Ng in 1987 and by Maksa and Páles in 2009, respectively.
\end{abstract}
	
\maketitle
	
\section{Introduction} 
The simplex of strictly ordered $n$-tuples of a set $H\subset \R$ denoted by $\sigma_n(H)$ is defined by
\Eq{*}{
\sigma_n(H) :=\{(x_1,\dots,x_n)\in H\quad |\quad x_1<\dots<x_n\}.
}
Obviously, the set $\sigma_n(H)$ is nonempty set if and only if the cardinality $|H|$ of $H$ is bigger than or equal to $n$. We adopt that $|H|\geq n$. Let $\omega=(\omega_1,\dots,\omega_n):H\to\R^n$ be a vector-valued function, and define the functional operator $\Phi_\omega:\sigma_n(H)\to \R$ by
\Eq{*}{
\Phi_\omega(x_1,\dots,x_n):=
\begin{vmatrix}
	\omega_1(x_1) \dots\omega_1(x_n)\\
	\vdots\quad\ddots\quad\vdots\\
	\omega_n(x_1)\dots \omega_n(x_n)
\end{vmatrix}
\qquad ((x_1,\dots,x_n)\in \sigma_n(H)).
}
A continuous function $\omega$ is said to be an \emph{$n$-dimensional positive (respectively negative) Chebyshev system over} $H$ if $\Phi_\omega$ is positive (respectively negative) over $\sigma_n(H)$. The system $\omega$ is called an \emph{$n$-dimensional Chebyshev system over} $H$ if it is either a positive or a negative Chebyshev system over $H$. If $\omega:\R\to\R^n$ equals the \emph{$n$-dimensional standard or polynomial system} $\pi_n:I\to\R^n$, which is defined by 
\Eq{*}{
 \pi_n(t):=(1,t,\dots,t^{n-1}) \qquad(t\in\R),
}
then, by computing Vandermonde determinants, one can easily show that it is a positive Chebyshev system. For standard results and applications related to Chebyshev systems, we refer to the monographs \cite{Kar68} and \cite{KarStu66}.

In what follows, we recall some definitions from the paper \cite{PalRad16} (see also the paper \cite{GilPal08} for these definitions in the polynomial setting). Let $I\subset \R$ be a nonvoid interval, $n\in\N$ and let $\omega:I\to\R$ be an $n$-dimensional positive Chebyshev system over $I$. 

For a given vector $t=(t_1,\dots,t_n)\in\R^n_+$, the function $f:I\to\R$ is said to be \emph{$(t,\omega)$-convex} if
\Eq{eq1}{
\Phi_{(\omega,f)}(x,x+t_1h,\dots,x+(t_1+\dots+t_n)h)\geq 0
}
holds for all $h>0$, $x\in I$ with $x+(t_1+\dots+t_n)h\in I$. If $T\subseteq\R_+$ and $f$ is $(t,\omega)$-convex for every $t\in T^n$, then $f$ is called \emph{$(T,\omega)$-convex}.

If $t=(t_1,\dots,t_n)\in\R^n_+$ and \eq{eq1} is satisfied with equality, then $f$ is called a \emph{$(t,\omega)$-affine function}. If $T\subseteq\R_+$ and $f$ is $(t,\omega)$-affine for every $t\in T^n$, then $f$ is called \emph{$(T,\omega)$-affine}. In particular, we say that $f$ is \emph{$\omega$-Jensen convex} if it is $(\{1\},\omega)$-convex, i.e., if
\Eq{oJ}{
\Phi_{(\omega,f)}(x,x+h,\dots,x+nh)\geq 0
}
holds for all $h>0$, $x\in I$ with $x+nh\in I$. If \eq{oJ} is valid with equality instead of inequality, then $f$ is said to be \emph{$\omega$-Jensen affine}. 

A function $f$ is termed \emph{$\omega$-convex} if it is $(\R_+,\omega)$-convex. It is easy to see that $f$ is $\omega$-convex on $I$ if and only if
\Eq{oc}{
\Phi_{(\omega,f)}(x_0,x_1,\dots,x_n)\geq 0 \qquad((x_0,x_1,\dots,x_n)\in\sigma_{n+1}(I)).
}
A function $f$ is called \emph{$\omega$-affine} if \eq{oc} is satisfied with equality.

We have to mention that in the case when $\omega=\pi_n$, then the concepts of $\omega$-convexity and $\omega$-Jensen convexity, was introduced by Hopf \cite{Hop26} and Popoviciu \cite{Pop44} (see also the book \cite{Kuc85} by Kuczma) and these properties were called convexity and Jensen convexity of order $(n-1)$, respectively.

The following result shows that, for any nonempty set $T\subseteq\R_+$, $(T,\omega)$-convexity implies $\omega$-Jensen convexity.

\Thm{t}{Let $T\subseteq\R_+$ be a nonempty set. If a function $f:I\to\R$ is $(T,\omega)$-convex (resp.\ $(T,\omega)$-affine), then it is $(\Q_+,\omega)$-convex (resp.\ $(\Q_+,\omega)$-affine), in particular, it is $\omega$-Jensen convex (resp.\ $\omega$-Jensen affine).
}
\begin{proof}
The result immediately follows from \cite[Theorem 5]{PalRad16}.
\end{proof}

In section $2$ we show that to any $\omega$-Jensen convex (resp. $\omega$-Jensen affine) function there exists a continuous $\omega$-convex (resp. $\omega$-affine) function so that these two functions coincide on dense subset of their domain. As a corollary, we obtain that $\omega$-convex functions are automatically continuous. The Bernstein-Doetsch theorem is generalized to the $\omega$-convexity with respect to Chebyshev system, that is, we show that an $\omega$-Jensen convex function which is bounded over some nonempty open subinterval, is also $\omega$-convex. We also establish characterizations of $\omega$-Jensen affine functions in several settings. A sufficient condition on $\omega$ that ensures the existence of discontinuous $\omega$-Jensen affine (resp. $\omega$-Jensen convex) function is given as well. Some of the results of this section extend that of the paper \cite{Mat08c} by Matkowski.

In Section 3, we generalize the concept of standard Wright convexity to $\omega$-Wright convexity with respect to a positive Chebyshev system, and establish its relationship with $\omega$-convexity and $\omega$-Jensen convexity. We also generalize Ng's decomposition theorem to setting of $\omega$-Wright convexity with respect certain positive Chebyshev systems. Finally, in the two dimensional case, we show that $\omega$-Wright convexity could be equivalent to the Wright convexity for certain positive Chebyshev systems.

\section{Results on $\omega$-Jensen functions}

In what follows, if $D$ is a subset of $I$ and $f:D\to\R$, then $f$ is said to be \emph{locally uniformly continuous} if, for all compact subinterval $[a,b]$ of $I$ and for all $\varepsilon>0$, there exists $\delta>0$ such that, for all $x,y\in [a,b]\cap D$ with $|x-y|<\delta$, we have that $|f(x)-f(y)|<\varepsilon$.

\Lem{UC}{Let $D$ be a dense subset of $I$ and let $f:D\to\R$ be a locally uniformly continuous function. Then $f$ admits a continuous extension to $I$.}

The verification of this lemma is standard, however, for the convenience of the reader, we provide a proof.

\begin{proof}
Let $x\in I$ be a fixed and let $(x_n)$ be an arbitrary sequence in $D$ such that $x_n\to x$ as $n\to\infty$. We point out first that $(f(x_n))$ is a Cauchy sequence.
The set $\{x_n:n\in\N\}\cup\{x\}$ is compact, therefore there exists a compact subinterval $[a,b]\subset I$ which contains all the members of the sequence $(x_n)$. 

Let $\varepsilon>0$ be arbitrary. By the local uniform continuity of $f$, there exists $\delta>0$ such that $|f(u)-f(v)|<\varepsilon$ for all $u,v\in [a,b]\cap D$ with $|u-v|<\delta$. The sequence $(x_n)$ is a Cauchy sequence, there exists $N\in\N$ such that, for all $n,m\geq N$, we have $|x_n-x_m|<\delta$. Consequently, 
\Eq{Cau}{
|f(x_n)-f(x_m)|<\varepsilon
} whenever for all $n,m\geq N$. Therefore $(f(x_n))$ is a Cauchy sequence which implies its convergence. Denote its limit by $g(x)$.

We show that $g$ is well defined at $x$. If $(x_n)$ and $(y_n)$ are two sequences in $D$ converging to $x$, then the sequence $(x_1,y_1,x_2,y_2,x_3,y_3\dots)$ is also converges to $x$. What we have proved above implies that $(f(x_1),f(y_1),f(x_2),f(y_2),f(x_3),f(y_3),\dots)$ is a convergent sequence and hence the subsequences $(f(x_n)$ and $(f(y_n))$ must have the same limit. Therefore the value of $g$ at $x$ does not depend on the choice of the sequence converging to $x$.

Clearly, if $x\in D$, then the sequence $x_n=x$ converges to $x$, hence $g(x)=f(x)$. Therefore, $g$ is an extension of $f$. 

In order to prove that $g$ is continuous on $I$, we show that it is locally uniformly continuous on $I$. To see this, let $[a,b]\subseteq I$. Let $\varepsilon>0$ be arbitrary and, using the local uniform continuity of $f$, choose $\delta>0$ such that, for all $u,v\in[a,b]\cap D$ with $|u-v|<\delta$, we have $|f(u)-f(v)|\leq\varepsilon$. 

Let $x,y\in[a,b]$ be arbitrary with $|x-y|<\delta$ and we are going to show that $|g(x)-g(y)|\leq\varepsilon$. Indeed, consider two sequences $(x_n)$ and $(y_n)$ in $[a,b]\cap D$ such that $x_n\to x$ and $y_n\to y$ as $n\to\infty$. Then $|x_n-y_n|\to|x-y|<\delta$, hence, there exists $N\in\N$ such that, for all $n\geq N$, we have $|x_n-y_n|<\delta$.
Therefore, according to the choice of $\delta$ (with $u:=x_n$ and $v:=y_n)$, we get that $|f(x_n)-f(y_n)|<\varepsilon$ holds for all $n\geq N$. Upon taking the limit $n\to\infty$, it follows that $|g(x)-g(y)|\leq\varepsilon$.
\end{proof}

\Thm{T1}{Let $n\geq2$ and let $\omega=(\omega_1,\dots,\omega_n):I\to \R^n$ be an $n$-dimensional positive Chebyshev system and let $\K$ be a subfield of $\R$. If a function $f:I\to\R$ is $(\K_+,\omega)$-convex (resp.\ $(\K_+,\omega)$-affine), then there exists a continuous $\omega$-convex (resp.\ $\omega$-affine) function $g:I\to\R$ such that $g|_{I\cap\K}=f|_{I\cap\K}$.}

\begin{proof}
If $f$ is $(\K_+,\omega)$-convex, then by definition we have the inequality
\Eq{conv}{
\Phi_{(\omega,f)}(x_0,x_1,\dots,x_n)\geq 0 
\qquad ((x_0,\dots,x_n)\in \sigma_{n+1}(I\cap\K)).
}
Indeed, apply the inequality \eq{eq1} with $x:=x_0$, $h:=1$, and $t_i:=x_i-x_{i-1}\in\K_+$. If $f$ is assumed to be $(\K_+,\omega)$-affine, then \eq{conv} is satisfied with equality.

Using the continuity of the function $\omega$, we are going to show that the restricted function $f|_{I\cap\K}$ is locally uniformly continuous. To prove this let $[a,b]\subset I$ be arbitrary. Without loss of generality, we can assume that $a,b\in\K$. We fix some elements $u_0<u_1<\dots<u_{n-2}<a$ and $b<v$ of the set $I\cap\K$. Then, for $x,y\in[a,b]\cap\K$ with $x<y$, we have that $(u_0,\dots,u_{n-2},x,y),(u_0,\dots,u_{n-3},x,y,v)\in\sigma_{n+1}(I\cap\K)$, therefore \eq{conv} implies
\Eq{Ineqs}{
\Phi_{(\omega,f)}(u_0,\dots,u_{n-2},x,y)\geq 0 \qquad\mbox{and}\qquad
\Phi_{(\omega,f)}(u_0,\dots,u_{n-3},x,y,v)\geq0.
}
The first of these inequalities implies that
\Eq{*}{
\begin{vmatrix}
\omega_1(u_0)&\dots&\omega_1(u_{n-2})&\omega_1 (x) & \omega_1(y) \\ 
\vdots & \ddots & \vdots & \vdots & \vdots \\ 
\omega_n(u_0)&\dots&\omega_n(u_{n-2})&\omega_n (x) & \omega_n(y) \\ 
f(u_0)&\dots&f(u_{n-2})& f(x) & f(y) 
\end{vmatrix}\geq 0.
}
Developing this determinant by the last row, this inequality is equivalent to
\Eq{PR}{
  f(y)P(x)-f(x)P(y)\geq R(x,y) \qquad((x,y)\in \sigma_2(I\cap\K)),
}
where $P:[a,b]\to\R$ and $R:[a,b]^2\to\R$ are defined by
\Eq{*}{
  P(z):=\Phi_\omega(u_0,\dots,u_{n-2},z)
}
and 
\Eq{*}{
  R(x,y):=\sum_{i=0}^{n-2} (-1)^{n-1-i} f(u_i) 
      \Phi_\omega(u_0,\dots,u_{i-1},u_{i-1},\dots,u_{n-2},x,y)
}
For $z\in[a,b]$, we have that $(u_0,\dots,u_{n-2},z)\in\sigma_n(I)$. Therefore, by the positivity and continuity of the Chebyshev system, it follows that $P$ is positive and continuous over $[a,b]$. We can also see that $R$ is continuous on $[a,b]^2$ and $R(z,z)=0$ for all $z\in[a,b]$. Substituting $x=a$ into \eq{PR}, we get that
\Eq{*}{
  f(y)\geq \frac{R(a,y)+f(a)P(y)}{P(a)} \qquad(y\in[a,b]\cap\K).
}
The right hand side of this inequality is a continuous function of $y$ over the compact interval $[a,b]$, therefore it is bounded from below. Hence $f|_{[a,b]\cap\K}$ is also bounded from below. Putting $y=b$ into \eq{PR}, it follows that
\Eq{*}{
  \frac{f(b)P(x)-R(x,b)}{P(b)}\geq f(x) \qquad(x\in[a,b]\cap\K).
}
Arguing similarly as above, this inequality yields that $f|_{[a,b]\cap\K}$ is bounded from above hence there exists a positive number $K$ such that, for $x\in[a,b]\cap\K$, we have $|f(x)|\leq K$.

Now we consider the second inequality in \eq{Ineqs}. It can be rewritten as
\Eq{*}{
\begin{vmatrix}
\omega_1(u_0)&\dots&\omega_1(u_{n-3})&\omega_1 (x) & \omega_1(y) & \omega_1(v)\\ 
\vdots & \ddots & \vdots & \vdots & \vdots & \vdots \\ 
\omega_n(u_0)&\dots&\omega_n(u_{n-3})&\omega_n (x) & \omega_n(y) & \omega_n(v) \\ 
f(u_0)&\dots&f(u_{n-3})& f(x) & f(y) & f(v)
\end{vmatrix}\geq 0.
}
Developing this determinant by its last row, this inequality is equivalent to
\Eq{QS}{
  f(x)Q(y)-f(y)Q(x)\geq S(x,y) \qquad((x,y)\in \sigma_2(I\cap\K)),
}
where $Q:[a,b]\to\R$ and $S:[a,b]^2\to\R$ are defined by
\Eq{*}{
  Q(z):=\Phi_\omega(u_0,\dots,u_{n-3},z,v)
}
and
\Eq{*}{
  S(x,y):=
     \sum_{i=0}^{n-3} (-1)^{n-1-i} f(u_i)
     \Phi_\omega(u_0,\dots,u_{i-1},u_{i-1},\dots,u_{n-3},x,y,v)
     -f(v)\Phi_\omega(u_0,\dots,u_{n-3},x,y).
}
The inclusion $(u_0,\dots,u_{n-3},z,v)\in\sigma_n(I\cap\K)$, the positivity and continuity of the Chebyshev system yield that $Q$ is a positive and continuous function over $[a,b]$. We also have that $S$ is continuous over $[a,b]^2$ and $S(z,z)=0$ for all $z\in[a,b]$. 

For $x,y\in[a,b]\cap\K$ the inequalities \eq{PR} and \eq{QS} imply that
\Eq{AB}{
  f(y)-f(x)
  &\geq \frac{P(y)-P(x)}{P(x)}f(x)+\frac{R(x,y)}{P(x)}
  \geq -K\frac{|P(y)-P(x)|}{P(x)}+\frac{R(x,y)}{P(x)}=:A(x,y),\\
  f(y)-f(x)
  &\leq \frac{Q(y)-Q(x)}{Q(x)}f(x)+\frac{S(x,y)}{Q(x)}
  \leq \frac{|Q(y)-Q(x)|}{Q(x)}K+\frac{S(x,y)}{Q(x)}=:B(x,y).
}
The functions $A$ and $B$ defined on the right hand side of these inequalities are continuous over $[a,b]^2$ and hence they are uniformly continuous over $[a,b]^2$. Therefore, for all $\varepsilon>0$, there exists $\delta>0$ such that, for all $(x,y),(u,v)\in[a,b]^2$ with $\|(x,y)-(u,v)\|<\delta$,
we have
\Eq{*}{
  \max(|A(x,y)-A(u,v)|,|B(x,y)-B(u,v)|)<\varepsilon.
}
In particular, if $|x-y|<\delta$, then substituting $(u,v)=(x,x)$ and using that $A$ and $B$ vanish at diagonal points of the square $[a,b]^2$, the above inequality implies that
\Eq{*}{
  \max(|A(x,y)|,|B(x,y)|)<\varepsilon.
}
Therefore, in view of the inequalities in \eq{AB}, $x,y\in[a,b]\cap\K$ with $|x-y|<\delta$, we obtain that
\Eq{*}{
  |f(y)-f(x)|<\varepsilon.
}
This proves that $f|_{[a,b]\cap\K}$ is uniformly continuous and hence $f|_{I\cap\K}$ is locally uniformly continuous. If $f$ is $(\K_+,\omega)$-affine, then it is also $(\K_+,\omega)$-convex, therefore we have the same conclusion.

In view of \lem{UC} with the dense set $D=I\cap\K$, there exists a continuous function $g:I\to\R$ such that $g(x)|_{I\cap\K}=f(x)|_{I\cap\K}$.

Finally, we show that $g$ is $\omega$-convex (resp.\ $\omega$-affine). Let $(y_0,\dots,y_n)$ be an arbitrary element of $\sigma_{n+1}(I)$. Then, by the density of $\K$ in $I$, for each $j\in\{0,\dots,n\}$, there exists  a sequence $(x_{k,j})_{k\in\N}$ in $\sigma_{n+1}(I\cap\K)$ converging to $y_j$ as $k\to\infty$. Then, applying the $(\K_+,\omega)$-convexity (resp.\ the $(\K_+,\omega)$-affinity) of $f$, we have that \eq{conv} is valid, we obtain
\Eq{*}{
	\Phi_{(\omega,g)}(x_{k,0},x_{k,1},\dots,x_{k,n})
	&=\Phi_{(\omega,f)}(x_{k,0},x_{k,1},\dots,x_{k,n})\geq0\\ 
	(\mbox{resp. } \Phi_{(\omega,g)}(x_{k,0},x_{k,1},\dots,x_{k,n})
	&=\Phi_{(\omega,f)}(x_{k,0},x_{k,1},\dots,x_{k,n})=0).
}
By the continuity of $g$, the function $\Phi_{(\omega,g)}$ is continuous. Upon taking the limit $k\to\infty$, the above inequality (resp.\ equality) implies that
\Eq{*}{
	\Phi_{(\omega,g)}(y_0,\dots,y_n)\geq0 \qquad(\mbox{resp. } \Phi_{(\omega,g)}(y_0,\dots,y_n)=0).
}
Therefore, $g$ is $\omega$-convex (resp.\ $\omega$-affine).
\end{proof}

\Cor{L0}{Let $n\geq2$ and $\omega=(\omega_1,\dots,\omega_n):I\to \R^n$ be an $n$-dimensional positive Chebyshev system. If $f:I\to\R$ is $\omega$-convex (resp.\ $\omega$-affine), then it is continuous on $I$.}

\begin{proof}
The corollary follows by applying \thm{T1} with $\K:=\R$.
\end{proof}

\Cor{L1}{Let $n\geq2$ and $\omega=(\omega_1,\dots,\omega_n):I\to \R^n$ be an $n$-dimensional positive Chebyshev system. If $f:I\to\R$ is $\omega$-Jensen convex (resp.\ $\omega$-Jensen affine), then there exists a continuous $\omega$-convex (resp.\ $\omega$-affine) function $g:I\to\R$ such that $g|_{I\cap\Q}=f|_{I\cap\Q}$.}

\begin{proof} In view of \thm{t}, the $\omega$-Jensen convexity (resp.\ $\omega$-Jensen affinity) of $f$ implies that it is $(\Q_+,\omega)$-convex (resp.\ $(\Q_+,\omega)$-affine). Now the statement of the corollary follows by applying \thm{T1} with $\K:=\Q$.
\end{proof}

The following statement is the extension of the celebrated Bernstein--Doetsch theorem \cite{BerDoe15} to the setting of $\omega$-Jensen convexity.

\Thm{BD}{If $f:I\to\R$ is $\omega$-Jensen convex and  bounded on a nonempty open subset of $I$, then it is continuous on $I$.}

\begin{proof}
Let $U$ be a nonvoid open subinterval of $I$ such that $f$ is bounded on $U$ by $K\geq0$. In the first part of the proof, we are going to show that $f$ is locally bounded on $I$, i.e., for every $v\in I$, there is an open set $V\subseteq I$ containing $v$ such that $f$ is bounded on $V$.

Let $v\in I$ be arbitrary. If $v\in U$, then the statement holds with $V=U$. Therefore, we may assume that $v\not\in U$. Choose a closed interval $[a,b]\subset U$. Then either $v<a$ or $b<v$. We consider now the case when $v<a$.

We choose some rational numbers $a-v< r_1<\dots<r_n<b-v$. Then we have that $v<a<v+r_1<\dots<v+r_n<b$. One can construct a bounded neighborhood $W$ of $v$ such that $\overline{W}\subseteq I$ and, for all $x\in W$, we have $x<a<x+r_1<\dots<x+r_n<b$. Now the $\omega$-Jensen convexity of $f$ and \thm{t} yield that $f$ is $(\Q_+,\omega)$-convex, hence we can get that
\Eq{*}{
  \Phi_{\omega,f}(x,x+r_1,\dots,x+r_n)\geq0
}
for all $x\in W$. This inequality implies that
\Eq{*}{
	\begin{vmatrix}
		\omega_1(x)&\omega_1(x+r_1)&\dots&\omega_1(x+r_n)\\ 
		\vdots  & \vdots & \ddots \\ 
		\omega_n(x)&\omega_n(x+r_1)&\dots&\omega_n(x+r_n)\\ 
		f(x)&f(x+r_1)&\dots&f(x+r_n) 
	\end{vmatrix}\geq 0.
}
Developing the determinant by last row, we obtain
\Eq{I1}{
	(-1)^nf(x)&\Phi_\omega(x+r_1,\dots,x+r_n)\\&+\sum_{i=1}^n (-1)^{n+i} f(x+r_i)\Phi_\omega(x,x+r_1,\dots,x+r_{i-1},x+r_{i+1},\dots,x+r_n)\geq 0.
}
By the boundedness of $f$ on $U$, the inclusion $x+r_i\in[a,b]\subseteq U$ we obtain that $(-1)^{n+i} f(x+r_i)\leq K$. The continuity of $\omega$ and positivity of Chebyshev system yield that the function
\Eq{*}{
  x\mapsto \sum_{i=1}^n\frac{  \Phi_\omega(x,x+r_1,\dots,x+r_{i-1},x+r_{i+1},\dots,x+r_n)}{\Phi_\omega(x+r_1,\dots,x+r_n)}
}
is bounded from above by a positive number $L$ over the compact set $\overline{W}$. Therefore, the inequality \eq{I1} implies that, for all $x\in W$, 
\Eq{*}{
	(-1)^{n-1}f(x)&\leq \sum_{i=1}^n (-1)^{n+i} f(x+r_i)\frac{\Phi_\omega(x,x+r_1,\dots,x+r_{i-1},x+r_{i+1},\dots,x+r_n)}{\Phi_\omega(x+r_1,\dots,x+r_n)}\leq KL,
}
which implies that $(-1)^{n-1}f$ is bounded from above over $W$.

To prove that $(-1)^{n-1}f$ is bounded from below over a neighborhood $V\subseteq W$ of $v$, we additionally fix $v_0\in I$ such that $v_0<v$. Now choose rational numbers 
\Eq{*}{
1<\frac{a-v_0}{v-v_0}< r_2<r_3<\dots<r_n<\frac{b-v_0}{v-v_0}.
}
Then we have
\Eq{*}{
a<v_0+r_2(v-v_0)<v_0+r_3(v-v_0)<\dots<v_0+r_n(v-v_0)<b.
}
Now one can construct a neighborhood $V$ of $v$ such that $V\subseteq W$ and, for all $x\in V$, we have
\Eq{*}{
v_0+r_0(x-v_0)<v_0+r_1(x-v_0)<a<v_0+r_2(x-v_0)<\dots<v_0+r_n(x-v_0)<b,
}
where $r_0=0$ and $r_1=1$.
Again applying \thm{t}, we conclude that
\Eq{*}{
\Phi_{(\omega,f)}(v_0+r_0(x-v_0),v_0+r_1(x-v_0),v_0+r_2(x-v_0),\dots,v_0+r_n(x-v_0))\geq0,
}
that is,
\Eq{*}{
\Phi_{(\omega,f)}(v_0,x,v_0+r_2(x-v_0),\dots,v_0+r_n(x-v_0))\geq0,
}
for all $x\in V$. This inequality, for $x\in V$, is equivalent to
\Eq{*}{
	\begin{vmatrix}
		\omega_1(v_0)&\omega_1(x)&\omega_1(v_0+r_2(x-v_0))&\dots&\omega_1(v_0+r_n(x-v_0))\\ 
		\vdots   & \vdots & \vdots & \ddots & \vdots \\ 
		\omega_n(v_0)&\omega_n(x)&\omega_n(v_0+r_2(x-v_0))&\dots&\omega_n(v_0+r_n(x-v_0))\\
		f(v_0)&f(x)&f(v_0+r_2(x-v_0))&\dots&f(v_0+r_n(x-v_0))\\ 
	\end{vmatrix}\geq 0.
}
For brevity, we write $s_i(x):=v_0+r_i(x-v_0)$, $i\in\{2,\dots,n\}$. Developing the determinant by its last row, we obtain
\Eq{*}{
(-1)^nf(v_0)&\Phi_\omega(x,s_2(x),\dots,s_n(x))+(-1)^{n-1}f(x)\Phi_\omega(v_0,s_2(x),\dots,s_n(x))\\
&+\sum_{i=2}^{n}(-1)^{n-i}f(s_i)\Phi_\omega(v_0,x,s_2(x)\dots,s_{i-1}(x),\dots,s_{i+1}(x),\dots,s_n(x))\geq0,
}
which yields
\Eq{I2}{
	(-1)^{n-1}f(x)&\geq (-1)^{n-1}f(v_0)\frac{\Phi_\omega(x,s_2(x),\dots,s_n(x))}{\Phi_\omega(v_0,s_2(x),\dots,s_n(x))}\\
	&\qquad+\sum_{i=2}^{n}(-1)^{n-i-1}f(s_i)\frac{\Phi_\omega(v_0,x,s_2(x)\dots,s_{i-1}(x),\dots,s_{i+1}(x),\dots,s_n(x))}{\Phi_\omega(v_0,s_2(x),\dots,s_n(x))}.
}
By the boundedness of $f$ on $U$ and the inclusion $s_i(x)\in[a,b]\subseteq U$, we get $(-1)^{n-i-1} f(s_i(x))\geq -K$. The continuity and the positivity of Chebyshev system $\omega$ yield that the functions
\Eq{*}{
	x\mapsto \frac{\Phi_\omega(x,s_2(x),\dots,s_n(x))}{\Phi_\omega(v_0,s_2(x),\dots,s_n(x))}\quad\mbox{and}\quad x\mapsto \sum_{i=1}^n\frac{  \Phi_\omega(v_0,x,s_2(x)\dots,s_{i-1}(x),\dots,s_{i+1}(x),\dots,s_n(x))}{\Phi_\omega(v_0,s_2(x),\dots,s_n(x))}
}
are bounded from above by positive numbers $M$ and $N$, respectively, over the compact set $\overline{V}\subseteq \overline{W}$. Therefore, the inequality \eq{I2} implies, for all $x\in V$, that
\Eq{*}{
	(-1)^{n-1}f(x)&\geq -|f(v_0)|M-KN,
}
which proves that $(-1)^{n-1}f$ is bounded from below over $V$. From the two-sided boundedness, it follows that $(-1)^{n-1}f$ is bounded over $V$, consequently, $f$ is also bounded over $V$.

To complete the proof of the theorem, we have to verify the continuity of $f$ at any point of $I$. Let $v\in I$ be arbitrary. Then, according to what we have proved in the first part, there exists a neighborhood $V\subseteq I$ of $v$ such that $f$ is bounded on $V$ by $K$. We are going to show that $f$ is uniformly continuous on every compact subinterval $[a,b]$ of $V$. This, in particular, implies the continuity of $f$ at $v$. 

Let $[a,b]\subseteq V$. We fix additional elements $a_1<b_1<\dots<a_{n-1}<b_{n-1}<a<b<a_{n}<b_{n}$ in $V$.
Then define the function $\Psi_1:[a_1,b_1]\times\dots\times[a_{n-1},b_{n-1}]\times[a,b]^2\to\R$ and $\Psi_2:[a_1,b_1]\times\dots\times[a_{n-2},b_{n-2}]\times[a,b]^2\times[a_n,b_n]\to\R$ by
\Eq{*}{
   \Psi_1(v_1,\dots,v_{n-1},x,y):=
   \bigg|1-\frac{\Phi_\omega(v_1,\dots,v_{n-1},y)}{\Phi_\omega(v_1,\dots,v_{n-1},x)}\bigg|
	+\sum_{i=1}^{n-1} \bigg|\frac{\Phi_\omega(v_1,\dots,v_{i-1},v_{i+1}\dots,v_{n-1},x,y)}{\Phi_\omega(v_1,\dots,v_{n-1},x)}\bigg|.
}
and
\Eq{*}{
	\Psi_2(v_1,\dots,v_{n-2},x,y,v_n):=
	&\bigg|1-\frac{\Phi_\omega(v_1,\dots,v_{n-2},x,v_n)}{\Phi_\omega(v_1,\dots,v_{n-2},y,v_n)}\bigg|+\frac{\Phi_\omega(v_1,\dots,v_{n-2},x,y)}{\Phi_\omega(v_1,\dots,v_{n-2},y,v_n)}\\
	&+\sum_{i=1}^{n-2} \bigg|\frac{\Phi_\omega(v_1,\dots,v_{i-1},v_{i+1}\dots,v_{n-2},x,y,v_n)}{\Phi_\omega(v_1,\dots,v_{n-2},y,v_n)}\bigg|.
}
Then the functions $\Psi_1$ and $\Psi_2$ are continuous over a compact rectangle, therefore, they are uniformly continuous. On the other hand, $\Psi_1(v_1,\dots,v_{n-1},x,x)=\Psi_2(v_1,\dots,v_{n-2},x,x,v_n)=0$ for all $x\in[a,b]$ and $(v_1,\dots,v_{n-1},v_n)\in [a_1,b_1]\times\dots\times[a_{n-1},b_{n-1}]\times[a_n,b_n]$.

Thus, for every $\varepsilon>0$, there exists $\delta>0$ such that, for all $x,y\in[a,b]$ with $|x-y|<\delta$ and for all $(v_1,\dots,v_{n-1},v_n)\in [a_1,b_1]\times\dots\times[a_{n-1},b_{n-1}]\times[a_n,b_n]$,
\Eq{*}{
  \Psi_1(v_1,\dots,v_{n-1},x,y)
  &=\Psi_1(v_1,\dots,v_{n-1},x,y)-\Psi_1(v_1,\dots,v_{n-1},x,x)
  <\frac{\varepsilon}{K},\\
	\Psi_2(v_1,\dots,v_{n-2},x,y,v_n)
	&=\Psi_2(v_1,\dots,v_{n-2},x,y,v_n)-\Psi_2(v_1,\dots,v_{n-2},x,x,v_n)
	<\frac{\varepsilon}{K}.
}

Now choose $x,y\in [a,b]$ with $x<y<x+\delta$. 
Then choose the rational numbers $\lambda_1,\dots,\lambda_n$ such that
\Eq{*}{
  \frac{a_i-x}{y-x}<\lambda_i<\frac{b_i-x}{y-x},\qquad i\in\{1,\dots,n\}.
}
Then, with the notation $v_i:=(1-\lambda_i)x+\lambda_i y$,
we have that
\Eq{*}{
  a_i<v_i<b_i, \quad i\in\{1,\dots,n\}.
}
First define 
\Eq{*}{
  r_i:=\frac{v_i-v_1}{y-v_1}\quad(i\in\{1,\dots,n-1\}), \qquad r_n:=\frac{x-v_1}{y-v_1},
  \qquad r_{n+1}:=\frac{y-v_1}{y-v_1}.
}
Clearly, due to the chain of inequalities $v_1<v_2\dots<v_{n-1}<x<y$, we have that $r_1=0<r_2<\dots<r_n<r_{n+1}=1$.
On the other hand, for $i\in\{1,\dots,n-1\}$, 
\Eq{*}{
  r_i=\frac{v_i-v_1}{y-v_1}
  =\frac{(1-\lambda_i)x+\lambda_i y-((1-\lambda_1)x+\lambda_1 y)}{y-((1-\lambda_1)x+\lambda_1 y))}
  =\frac{(\lambda_i-\lambda_1)(y-x)}{(1-\lambda_1)(y-x)}
  =\frac{\lambda_i-\lambda_1}{1-\lambda_1}\in\Q,
}
and
\Eq{*}{
  r_n=\frac{x-v_1}{y-v_1}
  =\frac{x-((1-\lambda_1)x+\lambda_1 y)}{y-((1-\lambda_1)x+\lambda_1 y))}
  =\frac{\lambda_1(x-y)}{(1-\lambda_1)(y-x)}
  =\frac{-\lambda_1}{1-\lambda_1}\in\Q.
}
Using \thm{t}, the $\omega$-Jensen convexity of $f$ implies that
\Eq{*}{
\Phi_{(\omega,f)}(v_1+r_1(y-v_1),\dots,v_1+r_{n-1}(y-v_1),v_1+r_n(y-v_1),v_1+r_{n+1}(y-v_1))\geq 0,
}
which, according to the definition of the numbers $r_0,r_1,\dots,r_{n-1},r_n$, yields that
\Eq{*}{
  \Phi_{(\omega,f)}(v_1,\dots,v_{n-1},x,y)\geq 0.
}
Now from the above inequality, we get
\Eq{*}{
	\begin{vmatrix}
		\omega_1(v_1)&\dots&\omega_1(v_{n-1})&\omega_1 (x) & \omega_1(y) \\ 
		\vdots & \ddots & \vdots & \vdots & \vdots \\ 
		\omega_n(v_1)&\dots&\omega_n(v_{n-1})&\omega_n (x) & \omega_n(y) \\ 
		f(v_1)&\dots&f(v_{n-1})& f(x) & f(y) 
	\end{vmatrix}\geq 0.
}
Developing this determinant by the last row, we obtain
\Eq{*}{
	&f(y)\Phi_\omega(v_1,\dots,v_{n-1},x)-f(x)\Phi_\omega(v_1,\dots,v_{n-1},y)\\
	&+\sum_{i=1}^{n-1} (-1)^{n-i-1}f(v_i)\Phi_\omega(v_1,\dots,v_{i-1},v_{i+1}\dots,v_{n-1},x,y)\geq 0.
}
By the positivity of Chebyshev system and the boundedness of $f$ over $V$, the above inequality yields
\Eq{xy}{
	f(x)-f(y)&\leq f(x)\bigg(1-\frac{\Phi_\omega(v_1,\dots,v_{n-1},y)}{\Phi_\omega(v_1,\dots,v_{n-1},x)}\bigg)\\
	&\qquad+\sum_{i=1}^{n-1} (-1)^{n-i-1}f(v_i)\frac{\Phi_\omega(v_1,\dots,v_{i-1},v_{i+1}\dots,v_{n-1},x,y)}{\Phi_\omega(v_1,\dots,v_{n-1},x)}\\
	&\leq K\Psi_1(v_1,\dots,v_{n-1},x,y)<\varepsilon.
}

Secondly define
\Eq{*}{
	&s_i:=\frac{v_i-v_1}{v_n-v_1}\quad(i\in\{1,\dots,n-2\}\cup\{n\}), \qquad s'_{n-1}:=\frac{x-v_1}{v_n-v_1}, \quad s''_{n-1}:=\frac{y-v_1}{v_n-v_1}.
}
Due to the chain inequalities $v_1<v_2\dots<v_{n-2}<x<y<v_n$, we have that $s_1=0<s_2<\dots<s'_{n-1}<s''_{n-1}<s_{n}=1$. On the other hand, for $i\in \{1,\dots,n-2\}\cup\{n\}$,
\Eq{*}{
s_i=\frac{v_i-v_1}{v_n-v_1}=\frac{(1-\lambda_i)x+\lambda_iy-((1-\lambda_1)x+\lambda_1y)}{(1-\lambda_n)x+\lambda_n y-((1-\lambda_1)x+\lambda_1y)}=\frac{(\lambda_i-\lambda_1)(y-x)}{(\lambda_n-\lambda_1)(y-x)}=\frac{\lambda_i-\lambda_1}{\lambda_n-\lambda_1}\in\Q,
}
\Eq{*}{
s'_{n-1}=\frac{x-v_1}{v_n-v_1}=\frac{x-((1-\lambda_1)x+\lambda_1y)}{(1-\lambda_n)x+\lambda_n y-((1-\lambda_1)x+\lambda_1y)}=\frac{-\lambda_1(y-x)}{(\lambda_n-\lambda_1)(y-x)}=\frac{-\lambda_1}{\lambda_n-\lambda_1}\in\Q
}
and
\Eq{*}{
s''_{n-1}=\frac{y-v_1}{v_n-v_1}=\frac{y-((1-\lambda_1)x+\lambda_1y)}{(1-\lambda_n)x+\lambda_n y-((1-\lambda_1)x+\lambda_1y)}=\frac{(1-\lambda_1)(y-x)}{(\lambda_n-\lambda_1)(y-x)}=\frac{1-\lambda_1}{\lambda_n-\lambda_1}\in\Q.
}
Now using \thm{t}, the $\omega$-Jensen convexity of $f$ implies that
\Eq{*}{
  \Phi_{(\omega,f)}(v_1+s_1(v_n-v_1),\dots,v_1+s'_{n-1}(v_n-v_1),v_1+s''_{n-1}(v_n-v_1),
    v_1+s_{n}(v_n-v_1))\geq 0,
}
which, according to the definition of the numbers $s_1,s_2,\dots,s'_{n-1},s''_{n-1},s_n$, yields that
\Eq{*}{
 \Phi_{(\omega,f)}(v_1,\dots,v_{n-2},x,y,v_n)\geq0.
}
This inequality can be written as
\Eq{*}{
	\begin{vmatrix}
		\omega_1(v_1)&\dots&\omega_1(v_{n-2})&\omega_1 (x) & \omega_1(y) & \omega_1(v_n)\\ 
		\vdots & \ddots & \vdots & \vdots & \vdots & \vdots \\ 
		\omega_n(v_1)&\dots&\omega_n(v_{n-2})&\omega_n (x) & \omega_n(y) & \omega_n(v_n) \\ 
		f(v_1)&\dots&f(v_{n-2})& f(x) & f(y) & f(v_n)
	\end{vmatrix}\geq 0.
}
Developing this determinant by the last row, we obtain
\Eq{*}{
	&f(v_n)\Phi_\omega(v_1,\dots,v_{n-2},x,y)-f(y)\Phi_\omega(v_1,\dots,v_{n-2},x,v_n)+f(x)\Phi_\omega(v_1,\dots,v_{n-2},y,v_n)\\
	&+\sum_{i=1}^{n-2} (-1)^{n-i-1}f(v_i)\Phi_\omega(v_1,\dots,v_{i-1},v_{i+1}\dots,v_{n-2},x,y,v_n)\geq 0,
}
by the positivity of Chebyshev system, the above inequality is equivalent to
\Eq{yx}{
	f(y)-f(x)&\leq f(y)\bigg(1-\frac{\Phi_\omega(v_1,\dots,v_{n-2},x,v_n)}{\Phi_\omega(v_1,\dots,v_{n-2},y,v_n)}\bigg)+f(v_n)\frac{\Phi_\omega(v_1,\dots,v_{n-2},x,y)}{\Phi_\omega(v_1,\dots,v_{n-2},y,v_n)}\\
	&\qquad+\sum_{i=1}^{n-2} (-1)^{n-i-1}f(v_i)\frac{\Phi_\omega(v_1,\dots,v_{i-1},v_{i+1}\dots,v_{n-2},x,y,v_n)}{\Phi_\omega(v_1,\dots,v_{n-2},y,v_n)}\\
	&\leq K\Psi_2(v_1,\dots,v_{n-2},x,y,v_n)<\epsilon.
}
Hence the inequalities \eq{xy} and \eq{yx} imply that, for all $x,y\in[a,b]$ with $|x-y|<\delta$,
\Eq{*}{
	|f(x)-f(y)|<\varepsilon
}
holds. This proves that $f$ is uniformly continuous on $[a,b]$. The closed interval $[a,b]\subseteq V$ was arbitrary, therefore $f$ is continuous on $V$.
\end{proof}

\Thm{CSol}{Let $f:I\to\R$ be a function which is bounded on a nonempty open subset of $I$. Then it is $\omega$-Jensen affine if and only if $f=\alpha_1\omega_1+\dots+\alpha_n\omega_n$ for some $\alpha_1,\dots,\alpha_n\in\R$.
}

\begin{proof}
Assume first that $f$ is an $\omega$-Jensen affine function. Then, it is $\omega$-Jensen convex, and using \thm{BD}, it follows that $f$ is continuous on $I$. We show first that, for all $x_0,x_1,\dots,x_n\in I$,
\Eq{0n}{
  \Phi_{\omega,f}(x_0,x_1,\dots,x_n)=0.
}
Indeed, if two of points $x_0,x_1,\dots,x_n$ coincide, then this equality is obvious. We may assume that these points are pairwise distinct moreover that $x_0<x_1<\dots<x_n$. Let $0<r_{1,k}<\dots<r_{n,k}$ be rational sequences converging to $0<x_1-x_0<\dots<x_n-x_0$, respectively. Then, according to \thm{t}, the $\omega$-Jensen affinity of $f$, for all $k\in\N$, yields that
\Eq{*}{
  \Phi_{\omega,f}(x_0,x_0+r_{1,k},\dots,x_0+r_{n,k})=0.
}
Using the continuity of $f$ and taking the limit $k\to\infty$, it follows that \eq{0n} holds.

Let us fix $x_1<\dots<x_n$ in $I$ arbitrarily. Then,
\Eq{*}{
  \Phi_{\omega,f}(x,x_1,\dots,x_n)=0.
}
holds for all $x\in I$, i.e.,
\Eq{*}{
	\begin{vmatrix}
		\omega_1(x)&\omega_1(x_1)&\dots&\omega_1(x_n)\\ 
		\vdots & \vdots & \ddots & \vdots \\ 
		\omega_n(x)&\omega_n(x_1)&\dots&\omega_n(x_n) \\ 
		f(x)&f(x_1)&\dots&f(x_n)
	\end{vmatrix}= 0.
}
Developing this determinant by the first column, we obtain
\Eq{*}{
\sum_{i=1}^n(-1)^{i-1}\omega_i(x)\Phi_{\omega_1,\dots,\omega_{i-1},\omega_{i+1},\dots,\omega_n,f}(x_1,\dots,x_n)+(-1)^n f(x)\Phi_\omega(x_1,\dots,x_n)=0.
}
Therefore,
\Eq{*}{
f(x)=\sum_{i=1}^n(-1)^{n-i}\frac{\Phi_{\omega_1,\dots,\omega_{i-1},\omega_{i+1},\dots,\omega_n,f}(x_1,\dots,x_n)}{\Phi_\omega(x_1,\dots,x_n)}\omega_i(x).
}
Now, with an obvious choice of $\alpha_1,\dots,\alpha_n\in\R$, we can see that $f=\alpha_1\omega_1+\dots+\alpha_n\omega_n$ holds.

The reversed statement is obvious, if $f$  is a linear combination of the coordinate functions of $\omega$, then $\Phi_{\omega,f}$ is identically zero by standard properties of determinants.
\end{proof}

The following question seems to be important: What is a necessary and sufficient condition on $\omega$ that ensures the existence of discontinuous $\omega$-Jensen affine or discontinuous $\omega$-Jensen convex functions? The following result provides a sufficient condition. To formulate and prove this sufficient condition we need the following lemma.

\Lem{Poly}{Assume that there exist a positive continuous function $\omega_0:I\to\R_+$ and a matrix $M=(a_{i,j})_{1\leq i\leq n,\,0\leq j\leq n-1}\in\R^{n\times n}$ such that
\Eq{omi}{
\omega_i(x)
 =(a_{i,n-1}x^{n-1}+\dots+a_{i,1}x+a_{i,0})\cdot\omega_0(x) \qquad(x\in I,\,i\in\{1,\dots,n\}).
}
Then $\det(M)>0$ and, for all $(x_1,\dots,x_n)\in\sigma_n(I)$,
\Eq{ot}{
  \Phi_\omega(x_1,\dots,x_n)
  =\omega_0(x_1)\cdots\omega_0(x_n)
    \cdot\det(M)\cdot \Phi_{\pi_n}(x_1,\dots,x_n).
}
Additionally, let $f:I\to\R$. Then, for all $(x_0,\dots,x_n)\in\sigma_{n+1}(I)$,
\Eq{oft}{
  \Phi_{\omega,f}(x_0,\dots,x_n)
  =\omega_0(x_0)\cdots\omega_0(x_n)\cdot\det(M)\cdot\Phi_{\pi_n,f/\omega_0}(x_0,\dots,x_n).
}}

\begin{proof}The equality in \eq{omi} and the product rule for determinants imply, for all $(x_1,\dots,x_n)\in\sigma_n(I)$, that
\Eq{*}{
  \Phi_\omega&(x_1,\dots,x_n)\\
  &=\begin{vmatrix}
	(a_{1,n-1}x_1^{n-1}+\dots+a_{1,1}x_1+a_{1,0})\omega_0(x_1)& \dots &(a_{1,n-1}x_n^{n-1}+\dots+a_{1,1}x_n+a_{1,0})\omega_0(x_n)\\
	\vdots & \ddots & \vdots \\
	(a_{n,n-1}x_1^{n-1}+\dots+a_{n,1}x_1+a_{n,0})\omega_0(x_1)& \dots &(a_{n,n-1}x_n^{n-1}+\dots+a_{n,1}x_n+a_{n,0})\omega_0(x_n)
	\end{vmatrix}\\
  &=\omega_0(x_1)\cdots\omega_0(x_n)\begin{vmatrix}
	a_{1,n-1}x_1^{n-1}+\dots+a_{1,1}x_1+a_{1,0}& \dots &a_{1,n-1}x_n^{n-1}+\dots+a_{1,1}x_n+a_{1,0}\\
	\vdots & \ddots & \vdots \\
	a_{n,n-1}x_1^{n-1}+\dots+a_{n,1}x_1+a_{n,0}& \dots &a_{n,n-1}x_n^{n-1}+\dots+a_{n,1}x_n+a_{n,0}
	\end{vmatrix}\\
 &=\omega_0(x_1)\cdots\omega_0(x_n)
    \begin{vmatrix}
    a_{1,0} & \dots & a_{1,n-1} \\
	\vdots & \ddots & \vdots \\
    a_{n,0} & \dots & a_{n,n-1}
   \end{vmatrix}\cdot
    \begin{vmatrix}
	x_1^0 & \dots &x_n^0\\
	\vdots & \ddots & \vdots \\
	x_1^{n-1} & \dots &x_n^{n-1}
	\end{vmatrix}\\
	&=\omega_0(x_1)\cdots\omega_0(x_n)\det(M)\cdot\Phi_{\pi_n}(x_1,\dots,x_n),
}
which proves \eq{ot}. 

The value of the determinant $\Phi_{\pi_n}(x_1,\dots,x_n)$ equals $\prod_{1\leq i<j\leq n}(x_j-x_i)>0$ (because it is of Vandermonde-type). Therefore, the positivity of the Chebyshev system $\omega$, the positivity of the function $\omega_0$ and the equality \eq{ot} yield that $\det(M)>0$. 

The equalities in \eq{omi} and the product rule for determinants again imply, for all $(x_0,\dots,x_n)\in\sigma_{n+1}(I)$, that
\Eq{*}{
  \Phi_{\omega,f}&(x_0,\dots,x_n)\\
  &=\begin{vmatrix}
	(a_{1,n-1}x_0^{n-1}+\dots+a_{1,1}x_0+a_{1,0})\omega_0(x_0)& \dots &(a_{1,n-1}x_n^{n-1}+\dots+a_{1,1}x_n+a_{1,0})\omega_0(x_n)\\
	\vdots & \ddots & \vdots \\
	(a_{n,n-1}x_0^{n-1}+\dots+a_{n,0}x_1+a_{n,0})\omega_0(x_0)& \dots &(a_{n,n-1}x_n^{n-1}+\dots+a_{n,1}x_n+a_{n,0})\omega_0(x_n)\\
	f(x_0) & \dots & f(x_n)
	\end{vmatrix}\\
  &=\omega_0(x_0)\cdots\omega_0(x_n)
    \begin{vmatrix}
	a_{1,n-1}x_0^{n-1}+\dots+a_{1,1}x_0+a_{1,0}& \dots &a_{1,n-1}x_n^{n-1}+\dots+a_{1,1}x_n+a_{1,0}\\
	\vdots & \ddots & \vdots \\
	a_{n,n-1}x_0^{n-1}+\dots+a_{n,1}x_0+a_{n,0}& \dots &a_{n,n-1}x_n^{n-1}+\dots+a_{n,1}x_n+a_{n,0}\\
	\frac{f}{\omega_0}(x_0) & \dots & \frac{f}{\omega_0}(x_n)
	\end{vmatrix}\\
	&=\omega_0(x_0)\cdots\omega_0(x_n)
	\begin{vmatrix}
    a_{1,0} & \dots & a_{1,n-1} & 0 \\
	\vdots & \ddots & \vdots \\
    a_{n,0} & \dots & a_{n,n-1} & 0 \\
    0 & \dots & 0 & 1 
   \end{vmatrix}\cdot
    \begin{vmatrix}
	x_0^0 & \dots & x_n^0\\[2mm]
	\vdots & \ddots & \vdots \\
	x_0^{n-1} & \dots & x_n^{n-1}\\[2mm]
	\frac{f}{\omega_0}(x_0)& \dots & \frac{f}{\omega_0}(x_n)
	\end{vmatrix}\\
	&=\omega_0(x_0)\cdots\omega_0(x_n)\cdot\det(M)\cdot\Phi_{\pi_n,f/\omega_0}(x_0,\dots,x_n),
}
which shows the validity of \eq{oft}.
\end{proof}

\Thm{GS}{Assume that there exist a positive continuous function $\omega_0:I\to\R_+$ and a matrix $M=(a_{i,j})_{1\leq i\leq n,\,0\leq j\leq n-1}\in\R^{n\times n}$ such that \eq{omi} holds.
Then $f:I\to\R$ is an $\omega$-Jensen convex (resp.\ $\omega$-Jensen affine) function if and only if $\frac{f}{\omega_0}$ is a $\pi_n$-Jensen convex (resp.\ $\pi_n$-Jensen affine) function.}

\begin{proof} 
According to formula \eq{oft} of \lem{Poly}, for all $h>0$ and $x\in I\cap(I-nh)$, we have that
\Eq{*}{
  \Phi_{\omega,f}(x,x+h,\dots,x+nh)
  =\prod_{i=0}^n\omega_0(x+ih)\cdot\det(M)\cdot\Phi_{\pi_n,f/\omega_0}(x,x+h,\dots,x+nh).
}
Due to the positivity of $\det(M)$ and the positivity of the function $\omega_0$, it follows that the inequality
$\Phi_{\omega,f}(x,x+h,\dots,x+nh)\geq0$ holds if and only if $\Phi_{\pi_n,f/\omega_0}(x,x+h,\dots,x+nh)\geq0$ is valid. This shows that $f$ is $\omega$-Jensen convex if and only if $\frac{f}{\omega_0}$ is $\pi_n$-Jensen convex.

Similarly, $\Phi_{\omega,f}(x,x+h,\dots,x+nh)=0$ if and only if $\Phi_{\pi_n,f/\omega_0}(x,x+h,\dots,x+nh)=0$, which proves that $f$ is $\omega$-Jensen affine if and only if $\frac{f}{\omega_0}$ is $\pi_n$-Jensen affine.
\end{proof}

We need to recall the following characterization of $\pi_n$-Jensen affine functions.

\Thm{pJa}{A function $f:I\to\R$ is $\pi_n$-Jensen affine  if and only if there exist a constant $A_0\in\R$, an additive function $A_1:\R\to\R$, a symmetric biadditive function $A_2:\R^2\to\R$, \dots, and a symmetric $(n-1)$-additive function $A_{n-1}:\R^{n-1}\to\R$ such that
\Eq{gensol0}{
  f(x)=A_{n-1}(x,\dots,x)+\dots+A_2(x,x)+A_1(x)+A_0 \qquad(x\in I).
}}

\Cor{GS}{Assume that there exist a positive continuous function $\omega_0:I\to\R_+$ and a matrix $(a_{i,j})_{1\leq i\leq n,\,0\leq j\leq n-1}\in\R^{n\times n}$ such that
\eq{omi} holds. Then $f:I\to\R$ is an $\omega$-Jensen affine function if and only if there exist a constant $A_0\in\R$, an additive function $A_1:\R\to\R$, a symmetric biadditive function $A_2:\R^2\to\R$, \dots, and a symmetric $(n-1)$-additive function $A_{n-1}:\R^{n-1}\to\R$ such that
\Eq{gensol}{
  f(x)=\big(A_{n-1}(x,\dots,x)+\dots+A_2(x,x)+A_1(x)+A_0\big)\omega_0(x) \qquad(x\in I).
}}

\begin{proof} Assume first that $f$ is $\omega$-Jensen affine. Then, by \thm{GS}, $\frac{f}{\omega_0}$ is $\pi_n$-Jensen affine. Hence, according to \thm{pJa}, there exist a constant $A_0\in\R$, and additive function $A_1:\R\to\R$, a symmetric biadditive function $A_2:\R^2\to\R$, \dots, a symmetric $(n-1)$-additive function $A_{n-1}:\R^{n-1}\to\R$ such that
\Eq{*}{
  \frac{f}{\omega_0}(x)=A_{n-1}(x,\dots,x)+\dots+A_2(x,x)+A_1(x)+A_0 \qquad(x\in I).
}
This proves that $f$ is of the form \eq{gensol}.

To prove the reversed implication, assume that there exist a constant $A_0\in\R$, and additive function $A_1:\R\to\R$, a symmetric biadditive function $A_2:\R^2\to\R$, \dots, a symmetric $(n-1)$-additive function $A_{n-1}:\R^{n-1}\to\R$ such that \eq{gensol} holds. Then, according to \thm{pJa}, 
$\frac{f}{\omega_0}$ is a $\pi_n$-Jensen affine function.
In view of \thm{GS}, this implies that $f$ is $\omega$-Jensen affine.
\end{proof}

\section{Wright convexity with respect to extended Chebyshev systems}

In 1954, Wright \cite{Wri54} introduced a concept of convexity which is stronger than Jensen convexity and weaker than convexity. A function $f:I\to\R$ is called \emph{Wright convex} if
\Eq{*}{
  f(tx+(1-t)y)+f((1-t)x+ty)\leq f(x)+f(y)
  \qquad(x,y\in I,\, t\in[0,1]).
}
A characterization and the ultimate understanding of Wright convexity was established by Ng \cite{Ng87b}, who proved that $f:I\to\R$ is Wright convex if and only if it is of the form $f=g+A|_I$, where $g:I\to\R$ is convex and $A:\R\to\R$ is additive. 

The higher-order generalizations of Wright convexity were introduced by Gilányi and Páles \cite{GilPal08} as follows. A function $f:I\to\R$ was called \emph{Wright convex of order $(n-1)$}, if 
\Eq{*}{
  \Delta_{h_1}\cdots\Delta_{h_n} f(x)\geq0
}
holds for all $h_1,\dots,h_n>0$ and $x\in I\cap(I-(h_1+\dots+h_n))$. One can easily see that Wright convex of order $1$ is equivalent to Wright convexity in the standard sense.

In what follows, we extend the notion of higher-order Wright convexity to the setting of Chebyshev systems. 
Let $\omega=(\omega_1,\dots,\omega_n):I\to\R^n$ be a positive $n$-dimensional Chebyshev system. We say that $\overline{\omega}:I\to\R^{n+1}$ \emph{is an extension of $\omega$} if there exists a continuous function $\omega_{n+1}:I\to\R$ such that $\overline{\omega}:=(\omega_1,\dots,\omega_n,\omega_{n+1})$ and $\overline{\omega}$ is a positive $(n+1)$-dimensional Chebyshev system. 

Let $\omega$ be a positive $n$-dimensional Chebyshev system
and $\overline{\omega}$ be an arbitrarily fixed extension of $\omega$. We say that a function \emph{$f:I\to\R$ is $\overline{\omega}$-Wright convex} if, for all $h_1,\dots,h_n>0$ and $x\in I\cap(I-(h_1+\dots+h_n))$, the inequality
\Eq{Wf}{
  \sum_{(i_1,\dots,i_n)}
  \frac{\Phi_{(\omega,f)}(x,x+h_{i_1},\dots,x+h_{i_1}+\dots+h_{i_n})}{\Phi_{\overline{\omega}}(x,x+h_{i_1},\dots,x+h_{i_1}+\dots+h_{i_n})}\geq0
}
holds, where the summation is taken over all permutation $(i_1,\dots,i_n)$ of the elements $\{1,\dots,n\}$.

Our first result establishes the connections between 
$\omega$-convexity, $\overline{\omega}$-Wright convexity and $\omega$-Jensen convexity.

\Thm{JWC}{Let $\omega$ be a positive $n$-dimensional Chebyshev system and $\overline{\omega}$ be an extension of $\omega$. Then every $\omega$-convex function is $\overline{\omega}$-Wright convex and every $\overline{\omega}$-Wright convex function is $\omega$-Jensen convex.}

\begin{proof} Assume first that $f:I\to\R$ is $\omega$-convex.
Then, for all $h_1,\dots,h_n>0$ and $x\in I\cap(I-(h_1+\dots+h_n))$ all permutation $(i_1,\dots,i_n)$ of the elements $\{1,\dots,n\}$, we have that
\Eq{*}{
  \Phi_{(\omega,f)}(x,x+h_{i_1},\dots,x+h_{i_1}+\dots+h_{i_n})\geq0.
}
On the other hand, by the positivity of the Chebyshev system $\overline{\omega}$, we also have that
\Eq{*}{
  \Phi_{\overline{\omega}}(x,x+h_{i_1},\dots,x+h_{i_1}+\dots+h_{i_n})>0.
}
These inequalities yield that \eq{Wf} is valid on the domain indicated and hence $f$ is $\overline{\omega}$-Wright convex.

To verify the second assertion, assume that $f:I\to\R$ is $\overline{\omega}$-Wright convex. Taking $h_1:=\dots=h_n:=h>0$ in inequality \eq{Wf}, for all $h>0$ and $x\in I\cap(I-nh)$, we obtain that
\Eq{Jf}{
  \frac{\Phi_{(\omega,f)}(x,x+h,\dots,x+nh)}{\Phi_{\overline{\omega}}(x,x+h,\dots,x+nh)}\geq0.
}
Due to the positivity of the Chebyshev system $\overline{\omega}$, it follows that  
$\Phi_{\overline{\omega}}(x,x+h,\dots,x+nh)$ is positive, therefore, we can conclude that $\Phi_{(\omega,f)}(x,x+h,\dots,x+nh)\geq0$, which shows that $f$ is $\omega$-Jensen convex.
\end{proof}

The next theorem describes the connection between $\overline{\omega}$-Wright convexity and Wright convexity of order $(n-1)$.

\Thm{ChWC}{Assume that there exist a positive continuous function $\omega_0:I\to\R_+$ and a matrix $M:=(a_{i,j})_{1\leq i\leq n,\,0\leq j\leq n-1}\in\R^{n\times n}$ such that \eq{omi} holds. Define $\omega_{n+1}:I\to\R$ by $\omega_{n+1}(t):=t^n\omega_0(t)$.
Then $\overline{\omega}:=(\omega,\omega_{n+1})$ is an extension of the Chebyshev system $\omega$. In addition, we have the following assertions:
\begin{enumerate}[(i)]
 \item For all $(x_0,x_1,\dots,x_n)\in\sigma_{n+1}(I)$, the equality
 \Eq{phq}{
  \frac{\Phi_{(\omega,f)}(x_0,x_1,\dots,x_n)}{\Phi_{\overline{\omega}}(x_0,x_1,\dots,x_n)}
  =[x_0,x_1,\dots,x_n;f/\omega_0]
 }
 holds. Furthermore, a function $f:I\to\R$ is $\omega$-convex if and only if $f/\omega_0$ is convex of order $(n-1)$.
 \item For all $h_1,\dots,h_n>0$ and $x\in I\cap(I-(h_1+\dots+h_n))$, the equality
\Eq{php}{
\sum_{(i_1,\dots,i_n)}
\frac{\Phi_{(\omega,f)}(x,x+h_{i_1},\dots,x+h_{i_1}+\dots+h_{i_n})}{\Phi_{\overline{\omega}}(x,x+h_{i_1},\dots,x+h_{i_1}+\dots+h_{i_n})}
= \frac{\Delta_{h_1}\cdots\Delta_{h_n} (f/\omega_0)(x)}{h_1\cdots h_n}
}
holds. Furthermore, a function $f:I\to\R$ is $\overline{\omega}$-Wright convex if and only if $f/\omega_0$ is Wright convex of order $(n-1)$.
\end{enumerate}}

\begin{proof} We first verify that $\overline{\omega}=(\omega,\omega_{n+1})$ is a positive Chebyshev system. Let $(x_0,x_1,\dots,x_n)\in\sigma_{n+1}(I)$. Then, applying the equality \eq{oft} of \lem{Poly} with $f:=\omega_{n+1}$, we get
\Eq{oot}{
  \Phi_{\overline{\omega}}(x_0,x_1,\dots,x_n)
  &=\omega_0(x_0)\cdots\omega_0(x_n)\cdot\det(M)
  \cdot \Phi_{(\pi_n,\omega_{n+1}/\omega_0)}(x_0,x_1,\dots,x_n)\\
  &=\omega_0(x_0)\cdots\omega_0(x_n)\cdot\det(M)
  \cdot \Phi_{\pi_{n+1}}(x_0,x_1,\dots,x_n)>0.
}
The last inequality is due to the fact that $\Phi_{\pi_{n+1}}(x_0,x_1,\dots,x_n)$ is a Vandermonde determinant and $x_0<x_1<\dots<x_n$. This proves that $\overline{\omega}$ is a positive Chebyshev system, indeed.

In the rest of the proof denote $g:=f/\omega_0$. 
To show that assertion (i) holds, let $(x_0,x_1,\dots,x_n)\in\sigma_{n+1}(I)$ be fixed. In view of the \lem{Poly}, we have the equality \eq{oft}. Combining this equality with \eq{oot}, we can obtain
\Eq{*}{
  \frac{\Phi_{(\omega,f)}(x_0,x_1,\dots,x_n)}{\Phi_{\overline{\omega}}(x_0,x_1,\dots,x_n)}
  =\frac{\Phi_{(\pi_n,g)}(x_0,x_1,\dots,x_n)}{\Phi_{\pi_{n+1}}(x_0,x_1,\dots,x_n)}.
}
From the theory of divided differences, we have the identity
\Eq{*}{
  \frac{\Phi_{(\pi_n,g)}(x_0,x_1,\dots,x_n)}{\Phi_{\pi_{n+1}}(x_0,x_1,\dots,x_n)}=[x_0,x_1,\dots,x_n;g],
}
which, together with the previous equality shows that \eq{phq} holds.

The function $f$ is $\overline{\omega}$-convex if and only if, for all $(x_0,x_1,\dots,x_n)\in\sigma_{n+1}(I)$, the left hand side of \eq{phq} is nonnegative. According to this equality, this happens if and only if the right hand side is nonnegative, i.e., if $f/\omega_0$ is convex of order $(n-1)$.

To show assertion (ii), let $h_1,\dots,h_n>0$ and $x\in I\cap(I-(h_1+\dots+h_n))$ be fixed. Therefore, with the substitutions $x_j:=x+h_{1}+\dots+h_{j}$, (where $j\in\{0,\dots,n\}$), \eq{phq} implies that
\Eq{*}{
  \frac{\Phi_{(\omega,f)}(x,x+h_{1},\dots,x+h_{1}+\dots+h_{n})}{\Phi_{\overline{\omega}}(x,x+h_{1},\dots,x+h_{1}+\dots+h_{n})}
  =[x,x+h_{1},\dots,x+h_{1}+\dots+h_{n};g].
}
Applying this equality for $(h_{i_1},\dots,h_{i_n})$ (instead of $(h_1,\dots,h_n)$), where $(i_1,\dots,i_n)$ is an arbitrary permutation of $(1,\dots,n)$, we can see that
\Eq{*}{
  \sum_{(i_1,\dots,i_n)}
  \frac{\Phi_{(\omega,f)}(x,x+h_{i_1},\dots,x+h_{i_1}+\dots+h_{i_n})}{\Phi_{\overline{\omega}}(x,x+h_{i_1},\dots,x+h_{i_1}+\dots+h_{i_n})}
  =\sum_{(i_1,\dots,i_n)} [x,x+h_{i_1},\dots,x+h_{i_1}+\dots+h_{i_n};g].
}
On the other hand, from the paper \cite{GilPal08}, we have that
\Eq{*}{
  \sum_{(i_1,\dots,i_n)} [x,x+h_{i_1},\dots,x+h_{i_1}+\dots+h_{i_n};g]=\frac{\Delta_{h_1}\dots\Delta_{h_n} g(x)}{h_1\cdots h_n}
}
holds, which, together with the previous equality implies \eq{php}.

The function $f$ is $\overline{\omega}$-Wright convex if and only if, for all $h_1,\dots,h_n>0$ and $x\in I\cap(I-(h_1+\dots+h_n))$, the left hand side of \eq{php} is nonnegative. According to this equality this happens to be valid if and only if the right hand side is nonnegative, i.e., if $f/\omega_0$ is Wright convex of order $(n-1)$.
\end{proof}

In the following result, we establish a characterization theorem for $\overline{\omega}$-Wright convexity provided that the underlying Chebyshev system is strongly related to the polynomial one. This result generalizes the decomposition theorem of Maksa and Páles \cite{MakPal09b} which is related to the polynomial system. An alternative and more elementary proof of that theorem has been recently given by the authors in \cite{PalShi22}.

\Thm{DWC}{Assume that there exist a positive continuous function $\omega_0:I\to\R_+$ and a matrix $M:=(a_{i,j})_{1\leq i\leq n,\,0\leq j\leq n-1}\in\R^{n\times n}$ such that \eq{omi} holds. Define $\omega_{n+1}:I\to\R$ by $\omega_{n+1}(t):=t^n\omega_0(t)$ and set $\overline{\omega}:=(\omega,\omega_{n+1})$. Then a function $f:I\to\R$ is $\overline{\omega}$-Wright convex if and only if there exist
an $\omega$-convex function $F:I\to\R$ and, for each $k\in\{1,\dots,n-1\}$, a symmetric $k$-additive mapping $A_k:\R^k\to\R$ and a real constant $A_0$ such that, for all $x\in I$,
\Eq{COM}{
  f(x)=F(x)+(A_0+A_1(x)+\dots+A_{n-1}(x,\dots,x))\omega_0(x).
}}
\begin{proof}
To prove the necessity, assume that the function $f$ is $\overline{\omega}$-Wright convex. Then the assertion (ii) of \thm{ChWC} implies that $f/\omega_0$ is Wright convex of order $(n-1)$. The decomposition theorem of higher order Wright convex functions \cite{MakPal09b} implies that
there exist a function $G:I\to\R$ which is convex of order $(n-1)$, a real constant $A_0$  and, for each $k\in\{1,\dots,n-1\}$, a symmetric $k$-additive mapping $A_k:\R^k\to\R$ such that, for all $x\in I$,
\Eq{foo}{
\frac{f}{\omega_0}(x)=G(x)+A_0+A_1(x)+\dots+A_{n-1}(x,\dots,x).
}
This implies that \eq{COM} holds with $F:=G\omega_0$ and $F/\omega_0$ is convex of order $(n-1)$. By assertion (i) of \thm{ChWC}, it follows that the function $F$ is $\omega$-convex.

To prove the sufficiency,  assume that \eq{COM} holds,  multiplying \eq{COM} by $1/\omega_0(x)$ implies that
\Eq{*}{
\frac{f}{\omega_0}(x)=\frac{F}{\omega_0}(x)+(A_0+A_1(x)+\dots+A_{n-1}(x,\dots,x)).
}
Since the function $F$ is $\omega$-convex, therefore the assertion (i) of \thm{ChWC} implies that $F/\omega_0$ is convex of order $(n-1)$. Again, by the decomposition theorem of higher order Wright convex functions \cite{MakPal09b}, we can conclude that $f/\omega_0$ is Wright convex of order $(n-1)$. Thus, the assertion (ii) of \thm{ChWC} yields that the function $f$ is $\overline{\omega}$-Wright convex.
\end{proof}

In our subsequent result we will prove that if the extension of the two dimensional polynomial system is not a polynomial of at most second degree, then the convexity with respect to the two dimensional polynomial system (i.e., standard convexity) is equivalent to Wright convexity with respect to this extension. For the proof of this result, we will need the following characterization of a polynomial of at most second degree.

\Lem{QP}{Let $\rho:I\to\R$ be a continuous function which satisfies the functional equation
\Eq{QP}{
  \frac{\rho(z)-\rho(y)}{z-y}
  =\frac{\rho(z+u)-\rho(y-u)}{(z+u)-(y-u)},
  \qquad u\geq 0,\, y,z\in (I+u)\cap(I-u),\, z+u\geq y.
}
Then $\rho$ is a polynomial of at most second degree over $I$.}

\begin{proof} If $u>0$ and $y\in(I+u)\cap(I-u)$, then the limit of the right hand side exists as $z\to y$ by the continuity of $\rho$, which shows that $\rho$ is differentiable at $y$ and we get 
\Eq{o3}{
 \rho'(y)=\frac{\rho(y+u)-\rho(y-u)}{2u},
 \qquad u>0,\, y\in (I+u)\cap(I-u).
}
Since $u>0$ was arbitrary, it follows that $\rho$ is differentiable everywhere on $I$. Now the right hand side of the above equality is differentiable with respect to $y$, which implies that $\rho$ is twice differentiable. Repeating this argument, it follows that $\rho$ is three times differentiable on $I$. We are going to show that $\rho'''$ is identically zero on $I$. 

Let $y\in I$ be fixed arbitrarily.
Rearranging the equation \eq{o3}, we can obtain that
\Eq{*}{
 2u\rho'(y)=\rho(y+u)-\rho(y-u),
 \qquad u\in (y-I)\cap(I-y).
}
Differentiating this equality three times with respect to $u$, we get that
\Eq{*}{
 0=\rho'''(y+u)+\rho'''(y-u),
 \qquad u>0,\, y\in (I+u)\cap(I-u).
}
With the substitution $u=0$, we conclude that $\rho'''(y)=0$. Therefore, $\rho'''$ is identically zero on $I$. This yields that $\rho$ has to be a polynomial of at most second degree.
\end{proof}

\Thm{o3}{Assume that there exist a positive continuous function $\omega_0:I\to\R_+$ and a matrix $M:=(a_{i,j})_{1\leq i\leq 2,\,0\leq j\leq 1}\in\R^{2\times 2}$ such that \eq{omi} holds for $n=2$. 
Assume that $\omega_3:I\to\R$ is a continuous function such that $\overline{\omega}=(\omega_1,\omega_2,\omega_3)$ is an extension of $\omega=(\omega_1,\omega_2)$ and $\omega_3/\omega_0$ is not  a polynomial of at most second degree. Then every $\overline{\omega}$-Wright convex function is $\omega$-convex, i.e., $\overline{\omega}$-Wright convexity is equivalent to $\omega$-convexity.}
\begin{proof}
Under the conditions of the theorem, assume that $f:I\to\R$ is an $\overline{\omega}$-Wright convex function. That is, the inequality
\Eq{*}{
	\frac{\Phi_{(\omega,f)}(x,x+h_1,x+h_1+h_2)}{\Phi_{\overline{\omega}}(x,x+h_1,x+h_1+h_2)}+\frac{\Phi_{(\omega,f)}(x,x+h_2,x+h_1+h_2)}{\Phi_{\overline{\omega}}(x,x+h_2,x+h_1+h_2)}\geq0
}
holds for all $h_1,h_2>0$ and $x\in I\cap(I-(h_1+h_2))$. Using \lem{Poly}, we can see that this inequality  is equivalent to
\Eq{*}{
&\frac{\omega_0(x)\omega_0(x+h_1)\omega_0(x+h_1+h_2)\cdot \det (M)\cdot\Phi_{\pi_2,f/\omega_0}(x,x+h_1,x+h_1+h_2)}{\omega_0(x)\omega_0(x+h_1)\omega_0(x+h_1+h_2)\cdot \det (M)\cdot\Phi_{\pi_ 2,\omega_3/\omega_0}(x,x+h_1,x+h_1+h_2)}
\\
&+\frac{\omega_0(x)\omega_0(x+h_2)\omega_0(x+h_1+h_2)\cdot \det (M)\cdot\Phi_{\pi_2,f/\omega_0}(x,x+h_2,x+h_1+h_2)}{\omega_0(x)\omega_0(x+h_2)\omega_0(x+h_1+h_2)\cdot \det (M)\cdot\Phi_{\pi_2,\omega_3/\omega_0}(x,x+h_2,x+h_1+h_2)}\geq0,
}
which simplifies to
\Eq{pi2}{
\frac{\Phi_{\pi_2,f/\omega_0}(x,x+h_1,x+h_1+h_2)}{\Phi_{\pi_ 2,\omega_3/\omega_0}(x,x+h_1,x+h_1+h_2)}
+\frac{\Phi_{\pi_2,f/\omega_0}(x,x+h_2,x+h_1+h_2)}{\Phi_{\pi_2,\omega_3/\omega_0}(x,x+h_2,x+h_1+h_2)}\geq0.
}
This means that $g:=f/\omega_0$ is $\overline{\pi}_2$-Wright convex, where $\overline{\pi}_2(t):=(1,t,\rho(t))$ and $\rho:=\omega_3/\omega_0$. 
Observe that, for $\varphi\in\{g,\rho\}$ and $i\in\{1,2\}$, we have
\Eq{*}{
  \Phi_{(\pi_2,\varphi)}(x,x+h_i,x+h_1+h_2)
  &=\begin{vmatrix}
    1 & 1 & 1 \\ x & x+h_i & x+h_1+h_2 \\
    \varphi(x) & \varphi(x+h_i) & \varphi(x+h_1+h_2)
   \end{vmatrix}\\[1mm] 
  &=h_{3-i}\varphi(x)-(h_1+h_2)\varphi(x+h_i)+h_i\varphi(x+h_1+h_1).
}
Using this formula, the inequality \eq{pi2} now states that
\Eq{fo}{
	&\frac{h_2g(x)-(h_1+h_2)g(x+h_1)+h_1g(x+h_1+h_2)}{h_2\rho(x)-(h_1+h_2)\rho(x+h_1)+h_1\rho(x+h_1+h_2)}\\
	&\qquad +\frac{h_1g(x)-(h_1+h_2)g(x+h_2)+h_2g(x+h_1+h_2)}{h_1\rho(x)-(h_1+h_2)\rho(x+h_2)+h_2\rho(x+h_1+h_2)}\geq 0
}
holds for all $h_1,h_2>0$ and $x\in I\cap(I-(h_1+h_2))$. 

By our assumptions, $(\pi_2,\rho)$ is an extension of $\pi_2$ and $\rho$ is not a polynomial of at most second degree. Using \lem{QP}, it follows that $\rho$ cannot satisfy the functional equation \eq{QP}, which means that there exist $u\geq 0$, $y,z\in (I+u)\cap(I-u)$ with $z+u\geq y$ such that
\Eq{*}{
	\frac{\rho(z)-\rho(y)}{z-y}
	\neq\frac{\rho(z+u)-\rho(y-u)}{(z+u)-(y-u)}.
}
Let $x:=y-u$, $h:=z-y+2u$, $t:=z-y+u$. Then $z=x+t$, $y=x+h-t$ and $z+u=x+h$, therefore the above relations state that 
\Eq{ht}{
	h(\rho(x+t)-\rho(x+h-t))+(h-2t)(\rho(x+h)-\rho(x))\neq 0
}
for some $h>0$, $x\in I\cap (I-h)$ and $t\in (0,h)$.

Let $h>0$ and $x\in I\cap (I-h)$ be fixed such that, for some $t\in(0,h)$, \eq{ht} holds. Define $T\subseteq (0,h)$ to be the set of those values $t$ for which \eq{ht} is valid. Then the set $T$ is nonempty and, by the continuity of $\rho$, it is also open. Let $T_+$ and $T_-$ denote the (disjoint) subsets of those elements $t\in T$, for which the left hand side of \eq{ht} is positive and negative, respectively. Then at least one of these subsets is nonempty (and also open).

Since $g$ is $\overline{\pi}_2$-Wright convex, hence \thm{JWC} implies that $g$ is $\pi_2$-Jensen convex, i.e., it is Jensen convex in the standard sense.
According to Rodé's Theorem \cite{Rod78}, $g$ is the pointwise maximum of Jensen affine functions, i.e., for all $p\in I$, there exists an additive function $A_p:\R\to\R$ such that
\Eq{gb}{
	g(y)\geq A_p(y-p)+g(p) \qquad(p,y\in I).  
}
Substituting $y:=x+t$ and $y:=x+h-t$, for $t\in T$ and $p\in I$, we get
\Eq{*}{
	g(x+t)\geq A_p(x+t-p)+g(p) \qquad g(x+h-t)\geq A_p(x+h-t-p)+g(p).
}
Therefore, with $h_1:=t$ and $h_2:=h-t$, the inequality \eq{fo} yields that
\Eq{*}{
&\frac{(h-t)g(x)-h\big(A_p(x+t-p)+g(p)\big)+tg(x+h)}{(h-t)\rho(x)-h\rho(x+t)+t\rho(x+h)}\\
&\qquad +\frac{tg(x)-h\big(A_p(x+h-t-p)+g(p)\big)+(h-t)g(x+h)}{t\rho(x)-h\rho(x+h-t)+(h-t)\rho(x+h)}\geq0.
}
Using the additivity of $A_p$ and moving the terms containing $A_p(t)$ to the right hand side, this inequality is equivalent to
\Eq{Ab}{
	&\frac{(h-t)g(x)-h\big(A_p(x-p)+g(p)\big)+tg(x+h)}{(h-t)\rho(x)-h\rho(x+t)+t\rho(x+h)}\\&\qquad +\frac{tg(x)-h\big(A_p(x+h-p)+g(p)\big)+(h-t)g(x+h)}{t\rho(x)-h\rho(x+h-t)+(h-t)\rho(x+h)}\\
	&\geq \frac{h(\rho(x+t)-\rho(x+h-t))+(h-2t)(\rho(x+h)-\rho(x))}{((h-t)\rho(x)-h\rho(x+t)+t\rho(x+h))(t\rho(x)-h\rho(x+h-t)+(h-t)\rho(x+h))}hA_p(t).
}
This inequality shows that $A_p$ is bounded from above on $T_+$ and is bounded from below on $T_-$. Therefore, $A_p$ is bounded from above or from below on a nonempty open subset of $T$. In view well-known properties of additive functions, this implies that $A_p$ is continuous, i.e., there exists a real constant $a_p$ such that $A_p(x)=a_px$ holds for all $x\in\R$. Thus, by \eq{gb}, we can see that $g$ is the pointwise maximum of continuous affine functions, Therefore, $g$ must be a convex function (in the standard sense). From this it follows that $f=g\omega_0$ is $\omega$-convex.
\end{proof}

It seems to be an open problem whether an analogue of the previous theorem is valid for the 3- or higher-dimensional setting.

\begin{thebibliography}{10}

\bibitem{BerDoe15}
F.~Bernstein and G.~Doetsch, \emph{{Zur {T}heorie der konvexen {F}unktionen}},
  Math. Ann. \textbf{76} (1915), no.~4, 514–526. \MR{1511840}

\bibitem{GilPal08}
A.~Gilányi and Zs. Páles, \emph{{On convex functions of higher order}}, Math.
  Inequal. Appl. \textbf{11} (2008), no.~2, 271–282.

\bibitem{Hop26}
E.~Hopf, \emph{{{Ü}ber die {Z}usammenhänge zwischen gewissen höheren
  {D}ifferenzenquotienten reeller {F}unktionen einer reellen {V}ariablen und
  deren {D}ifferenzierbarkeitseigenschaften}}, Ph.D. thesis,
  Friedrich–Wilhelms–Universität Berlin, 1926.

\bibitem{Kar68}
S.~Karlin, \emph{{Total positivity. {V}ol. {I}}}, Stanford University Press,
  Stanford, California, 1968. \MR{37 \#5667}

\bibitem{KarStu66}
S.~Karlin and W.~J. Studden, \emph{{Tchebycheff systems: {W}ith applications in
  analysis and statistics}}, {Pure and Applied Mathematics, Vol. XV},
  Interscience Publishers John Wiley \& Sons, New York-London-Sydney, 1966.
  \MR{34 \#4757}

\bibitem{Kuc85}
M.~Kuczma, \emph{{An {I}ntroduction to the {T}heory of {F}unctional {E}quations
  and {I}nequalities}}, {Prace Naukowe Uniwersytetu Śląskiego w Katowicach},
  vol. 489, Państwowe Wydawnictwo Naukowe — Uniwersytet Śląski,
  Warszawa–Kraków–Katowice, 1985, 2nd edn. (ed. by A. Gilányi),
  Birkhäuser, Basel, 2009. \MR{0788497 (86i:39008), MR 2467621}

\bibitem{MakPal09b}
Gy. Maksa and Zs. Páles, \emph{{Decomposition of higher-order {W}right-convex
  functions}}, J. Math. Anal. Appl. \textbf{359} (2009), 439–443.
  \MR{2010d:26014}

\bibitem{Mat08c}
J.~Matkowski, \emph{{Generalized convex functions and a solution of a problem
  of {Z}s. {P}áles}}, Publ. Math. Debrecen \textbf{73} (2008), no.~3-4,
  421–460. \MR{2466385 (2009i:26048)}

\bibitem{Ng87b}
C.~T. Ng, \emph{{{F}unctions generating {S}chur-convex sums}}, {General
  Inequalities, 5 (Oberwolfach, 1986)} (W.~Walter, ed.), {International Series
  of Numerical Mathematics}, vol.~80, Birkhäuser, Basel–Boston, 1987,
  p.~433–438. \MR{91b:26018}

\bibitem{Pop44}
T.~Popoviciu, \emph{{Les Fonctions Convexes}}, {Actualités Scientifiques et
  Industrielles, No. 992}, Hermann et Cie, Paris, 1944. \MR{0018705}

\bibitem{PalShi22}
Zs. Páles and M.K. Shihab, \emph{{Decomposition of higher-order {W}right
  convex functions revisited}}, Results Math. \textbf{77} (2022), no.~2, Paper
  No. 73, 11. \MR{4379106}

\bibitem{PalRad16}
Zs. Páles and É. {Székelyné Radácsi}, \emph{{Characterizations of
  higher-order convexity properties with respect to {C}hebyshev systems}},
  Aequationes Math. \textbf{90} (2016), no.~1, 193–210. \MR{3471295}

\bibitem{Rod78}
G.~Rodé, \emph{{Eine abstrakte {V}ersion des {S}atzes von {H}ahn–{B}anach}},
  Arch. Math. (Basel) \textbf{31} (1978), 474–481.

\bibitem{Wri54}
E.~M. Wright, \emph{{An inequality for convex functions}}, Amer. Math. Monthly
  \textbf{61} (1954), 620–622.

\end{thebibliography}

\providecommand{\bysame}{\leavevmode\hbox to3em{\hrulefill}\thinspace}
\def\MR#1{}
\providecommand{\MRhref}[2]{%
  \href{http://www.ams.org/mathscinet-getitem?mr=#1}{#2}
}
\providecommand{\href}[2]{#2}


\end{document}